\documentclass[reqno,11pt]{amsart}
\usepackage{amscd,amssymb,verbatim,array}
\usepackage{hyperref}
\usepackage{amsmath}
\usepackage[active]{srcltx}
\usepackage{tikz-cd}
\usepackage[all]{xy}

\numberwithin{equation}{section} \theoremstyle{plain}

\newtheorem{theorem}{Theorem}[section]

\newtheorem{proposition}[theorem]{Proposition}

\theoremstyle{definition}

\theoremstyle{remark}

\numberwithin{equation}{section}

\newcommand{\Vol}{\operatorname{Vol}}
\newcommand{\Tor}{\operatorname{Tor}}
\newcommand{\vol}{\operatorname{vol}}

\newcommand{\ddet}{\operatorname{det}}

\newcommand{\grad}{\operatorname{grad}}

\newcommand{\ind}{\operatorname{ind}}

\newcommand{\Id}{\operatorname{Id}}



\begin{document}
\title {Virtually small  spectral package 
of a Riemannian manifold}

\author{Dan Burghelea, Yoonweon Lee}



\date{\today}

\begin{abstract}
For a Morse function $f$ on a closed orientable Riemannian manifold $(M,g)$ one introduces the {\it virtually small spectral package},  an  analytic object consisting of a finite number of analytic quantities derived from $(g,f)$ which, in principle, can be calculated.
One shows that they  determine  the  {\it Torsion } of the underlying space $M,$  a parallel to the result that  the dimensions of the spaces of harmonic forms  calculate the Euler-Poincar\'e characteristic of $M$  and extend  the Poincar\'e duality  between harmonic forms and between Betti numbers for a closed oriented Riemannian manifold .
\end{abstract}

\maketitle

\setcounter{tocdepth}{1}
\tableofcontents

\section {Introduction}

For a compact ANR $X,$   the integral homology $H_r(X;\mathbb bZ)$ is a finitely generated abelian group of a finite rank $\beta_r(X)$ whose set of finite order elements has a finite cardinality $\Tor_i(X).$
The following two numbers are remarkable topological invariants
$$\chi(X):=\sum (-1)^i \beta_i(X)\ \quad \text {and} \quad \  {\mathbb T}or(X):=\prod (\Tor_i(X))^{(-1)^i}.$$

If $M^n$ is an $n$-dimensional closed orientable (topological) manifold, the Poincar\'e duality implies $\beta_r(M)= \beta_{n-r}(M)$ and $\Tor_r(M)= \Tor_{n-1-r}(M)$  and therefore for $n$ odd $\chi(M^n)=0$ and for $n$ even $ {\mathbb T}or(M^n)=1.$

If $(M,g)$ is an orientable closed Riemannian manifold and ${\mathcal H}^r(M,g)\subset \Omega^r(M)$ denotes the space of harmonic forms of degree $r$ and $\mathcal H^r_0 \subset \mathcal H^r$ the subspace of  integral forms,
then $T_{r}(M,g) := \mathcal H^r(M)/ \mathcal H^r_0(M)$ is the compact torus  of dimension $\beta_r(M)$ \footnote {in view of the Hodge de-Rham theorem }
as a Riemannian manifold with a metric induced from the scalar product on $\Omega^r(M)$ provided by $g$. Let $V_r(M,g):= \Vol (T_r(M,g))$ and $\mathbb V(M,g) :=\prod (V_{r}(M,g))^{(-1)^r} \in  \mathbb R_{> 0}.$  This is a Riemannian invariant. In view of the Hodge-de Rham theorem $\chi(M)$ is also a Riemannian invariant.

If  $f:M\to \mathbb R$ is a smooth real-valued function, $f$ defines a deformation (parametrization) $(\Omega^\ast(M), d^\ast (t))$ of the de Rham complex $(\Omega^\ast (M), d^\ast)$ with $d^\ast (t):= e^{-t f} d^\ast  e^{tf}= d^\ast + t dh\wedge \cdot,$ and therefore the one parameter family of second order elliptic differential operators  $\Delta_{g,f}^q(t): \Omega^q(M)\to \Omega^q(M)$  with
 $$ \Delta_{g,f}^{q}(t):= d^{q-1}(t)\cdot  \delta_g^q(t)  +  \delta_g^{q+1}(t) \cdot d^q(t)=  \Delta_g^q + t (L^q_X+ \mathcal L^q_X) + t^2||X||^2, $$
  where
 \begin{enumerate}
 \item $\delta_g^q(t): = (-1) ^{n(q-1) +1}\star_g^{n-q+1} \cdot  \ e^{tf} 
  d^{n-q} \ e^{-tf} \cdot \star_g ^q,$
 \item $X=- grad_gf$, $L_X$ the Lie derivative in direction $X,$ $\mathcal L_X^q:= (-1)^{(n+1) q+1} \star^{n-q} \cdot L^{n-q}_X \cdot \star^q,$
\item $||X||: M\to \mathbb R_{\geq 0}$ the length of the vector $X(x)= -\grad_g f(x).$
 \end{enumerate}
  These operators remain self-adjoint, nonnegative elliptic differential operators  on $L_2(\Omega^q(M)),$ the $L_2$-completion of $\Omega^q(M)$ with $\Delta_{g,f}^q(t)$ a zero order perturbation of the standard Laplace-Beltrami operator $\Delta_g^q$ \footnote {$L_X$ and $\mathcal L_X$ are order one differential operators but $L_X + \mathcal L_X$ is of order zero as well as the multiplication by the smooth function $||X||^2$}.

When  both $g$ and $f$ are implicit from the context, one abbreviates  $\Delta^q_{g,f} (t), \Delta^q_g, \delta^q_g, \star^q_g$ to $\Delta^q (t), \Delta^q, \delta^q, \star^q$ for simplicity in writing.

 The operators $\Delta^q(t),$ referred below as Witten Laplacians,  provide a {\it holomorphic family of type A } of self-adjoint operators in the sense of Kato  cf \cite{TK}
 and therefore, in view of a theorem of Reillich-Kato, cf. Theorem 3.9  chapter 7 in \cite {TK},  one has:
 \begin{theorem} (Rellich - Kato)\label {T3.1}\
 There exist a collection of non-negative real-valued  functions  $\lambda^q_\alpha(t),$ unique up to permutation, and a collection of  norm one $q$-differential form-valued maps  $\omega^q_\alpha(t) \in \Omega^q(M) ,$ analytic in $t\in \mathbb R,$  indexed by  $\alpha \in {\mathcal A}^{q},$ ${\mathcal A}^{q}$ a countable set, each with holomorphic extension to a neighborhood of the real line $\mathbb R\subset \mathbb C$ \footnote {holomorphic  extension means extensions $\lambda^q(z)\in \mathbb C$ , $\omega^q(z)\in \Omega(M)\otimes C$ for $z$ in a neighborhood of $\mathbb R$ in
 $\mathbb C$ which for $t\in \mathbb R$ is a real number and $\omega^q(t)\in \Omega(M)\otimes 1$ }
  such that:
 \begin{enumerate}
 \item $\Delta^q(t)\omega^{q}_\alpha (t)= \lambda^{q}_\alpha (t)\omega^{q}_\alpha (t),$
\item for any $t$ the collections $\lambda^{q}_\alpha (t)$ represent all repeated  eigenvalues of $\Delta^q(t)$  and the collection $\omega^{q}_\alpha (t)$ form a complete orthonormal family  of associated eigenvectors for the operator $\Delta^q(t),$
\item exactly $\beta_q= \dim H^q(M;\mathbb R)$ eigenvalue functions $\lambda^q_\alpha(t)$ are identically zero and all others are strictly positive.
\end{enumerate}
\end{theorem}

These analytic maps $\lambda^q_\alpha (t)$ and $\omega^q_\alpha (t)$ are called branches, eigenvalue branch and eigenform branch respectively, with extensions to holomorphic maps in the neighborhood of $\mathbb R\subset \mathbb C.$
The maps $\omega^q_\alpha(t)$ have $\pm$ ambiguity and if the branch $\lambda^q_\alpha(t)$ has multiplicity $\geq 2$, i.e. $\lambda^q_\alpha(t)= \lambda^q(t)$ for a finite set of indices $\alpha$, then
the family of finite dimensional vector spaces spanned by the corresponding $\omega^q_\alpha (t)$  is unique and of course analytic in $t.$

 If $f$ is a Morse function with $c_q$ critical points of Morse index $q$, then in view of a result of  Witten \cite {W}, cf \cite {BFKM} Proposition 5.2 or \cite {BFK2} Theorem 2.8
 for details,
 for any $q$ exactly  $c_q$ eigenvalue branches of $\Delta^q(t)$ go exponentially fast to zero and
 all others go at least  linearly fast to $\infty.$ Moreover each  eigenvalue branch  which converges to zero corresponds to a critical point and its corresponding eigenform branch  concentrates to
this critical point.
 We index these analytic functions as $\lambda^{\ind x}_x(t)$ and $\omega^{\ind x}_x(t)$ and refer to the finite collection of branches $$\boxed {\{\lambda^{\ind x}_x(t), \omega^{\ind x}_x(t), x\in Cr(f)\}}$$ as  the {\it virtually small spectral package of $(M,g,f)$} and to  the finite set $$\boxed{\{\lambda^{\ind x}_x(0), \omega^{\ind x}_x(0), x\in Cr(f)\}}$$ the subset of the infinite set $\{\lambda^q_\alpha, \omega^q_\alpha\},$ the spectral package of $(M,g),$ as the {\it virtually small spectral package of $(M,g)$ determined by $f$.}


Note 
that for $t$ large enough, in view of the {\it spectral gap theorem}, cf Theorem \ref{T32} stated in Section \ref{S3}, the eigenvalues $\lambda^q_x(t) x\in Cr_q(f),$ exhaust the first  $c_q= \sharp Cr_q(f),$ possibly repeated,  smallest  eigenvalues of $\Delta^q(t),$  however this is not true for $t=0$ as the Example  in section \ref{S3} shows. This explains the name {\it virtually small} for the collection $\{\lambda^{\ind x}_x(0), \omega^{\ind x}_x(0), x\in Cr(f)\}.$

 We denote by $\Omega^\ast_{vs}(M)(t)$ the span of the eigenforms $\omega^{\ind x}_x(t)  x\in Cr(f)$  inside
 $\Omega^\ast (M)$. They generate a finite dimensional sub-complex $(\Omega^\ast_{vs} (M)(t), d^\ast (t))$  of $(\Omega^\ast (M), d^\ast (t))$  with the $q$-component of dimension $c_q,$
 which is  an analytic family of cochain complexes with cohomology of constant dimension $\beta_r(M).$
 Since  $f$ is a Morse function, for any $x\in Cr(f)$  the stable/unstable set  $W^{\pm}_x$ of the vector field $-grad_g f$  are submanifolds diffeomorphic to $\mathbb R^{n-\ind x}/ \mathbb R^{\ind x}.$  For any $x$ choose an orientation $\mathcal O_x$ for $W^-_x$
and for any $x,y\in Cr_q(f)$ consider the integral
\begin{equation}\label {E1.1}
 A^q(x,y)(t)  :=  \int_{W^-_y} e^{tf}\omega^q_x(t),
 \end{equation}
 which a priory might not be convergent but  when convergent for any $x,y\in Cr_q(f)$ provide the non-negative number
$$a^q(M,g,f)(t):=  |\ \det |||A^q(x,y)(t)|||\  | \geq 0, $$
where $|||A^q(x,y)(t)|||$ is a $c_q\times c_q$ matrix with $c_q=\sharp Cr_q(f).$
Changing of the orientation $\mathcal O_y$ changes the sign of the integrals $A^q(x,y)$ for all $x$ but leave  $a^q(M,g,f)$ unchanged when defined.
\begin{proposition}\label {P1.2}\

\begin{enumerate}
\item If the vector field $-grad_gf$ is Morse-Smale, then  the integral (\ref{E1.1}) is uniformly convergent and both $A^q(x,y)(t)$ as well as $a^q(M,g,f)(t)$ are analytic in $t,$ with  the last being independent on the choice of $\omega^q_x(t)$ and the orientations $\mathcal O_x.$
\item The analytic function $a^q(M,g,f)(t)$  is  non-negative  with at most finitely many zeros in any interval $[T,\infty).$
\end{enumerate}
\end{proposition}


 In particular $$a(M,g,f)(t):=\prod (a^q(M,g,f)(t))^{(-1)^{q}}$$ is a priory a nonnegative
 meromorphic function in $t$ with at most finitely many zeros and poles in any interval $[T,\infty).$

The main result of this Note is the following theorem:

 \begin {theorem}\label {T1.3}
Suppose that the vector field $-\grad_g f$ is Morse-Smale. Then the following holds true.
 \begin {enumerate}
 \item The meromorphic map $a(M,g,f)(t)$ is strictly positive and has no zeros  and no poles.
 \item  The virtually small spectral package determined by $f$ together with the numbers $a(M,g,f)$ and ${\mathbb V}(M,g)$ are
 all analytic invariants and determine the topological invariant  ${\mathbb T}or(M)$ by the formula:
$$\log {\mathbb T}or (M)= 1/2 \sum _q   (-1)^{q+1} q \left( \sum_{\alpha\in {\mathcal A}^q_{vs, +}} \log \lambda^q_\alpha(0) \right)  + \log a(M,g,f)  - \log {\mathbb V}(M,g),$$
\end{enumerate}
where $a(M,g,f):= a(M,g,f)(0)$ and ${\mathcal A}^{q}_{vs, +} \subset {\mathcal A}^{q}$ is the set of indices $\alpha$'s such that $\lambda_{\alpha}^{q}(t)$ belongs to the virtually spectral package with $\lambda_{\alpha}^{q}(t) > 0$.
\end{theorem}

{\bf Conjecture 1:}  {\it The statement remains true without the hypothesis that "$-grad _g f$ is Morse-Smale".}
\vskip .1in

{\bf Conjecture 2:} {\it Under the hypothesis that " $-\grad_gf$ is Morse-Smale " one has  $a^q(M,g,f)(t)\ne 0$. }


Note that
%
if Conjecture 2 holds true for $t=0$, then it can be shown that the Morse complex defined by $(g,f)$ can be canonically realized as a sub-complex of the de Rham complex equipped with the scalar product defined by the metric $g.$ Recall that
Hodge-de Rham theorem implies that the  complex $(H^r(M), 0)$ can be realized in this way as the sub-complex of harmonic forms.

As shown in Section \ref{S4}, for an oriented closed Riemannian manifold the Hodge star operator $$\star : \Omega^q(M) \to \Omega^{n-q}(M)$$ identifies the virtually small $q-$spectral package of $(M,g, f)$ to the  the virtually small $(n-q)-$spectral package of $ (M,g, -f).$  This can be viewed as an extension of Poincar\'e duality.

\section {Proof of Proposition \ref{P1.2} and Theorem \ref{T1.3}}

Proof of Proposition \ref{P1.2}:

 One  says that the vector field $X=-grad_gf$  is {\it Morse-Smale} if for any $x,y\in Cr(f)$  the unstable set $W^-_x$ and the stable set  $W^+_y$ are transversal, which implies that $\mathcal T(x,y)= (W^-_x\cap W^+_y)/ \mathbb R$ \footnote {$\mathbb R$ acts freely by translation along the flow defined by $-\grad_gf$},
 the space of trajectories from $x$ to $y,$ is a manifold of dimension $\ind (x)- \ind (y)-1.$
Under the hypothesis that $X$ is Morse-Smale, it is shown in \cite{BFK} or  \cite{BH}
that the embedding $i_x: W^-_x\to M$ extends to a smooth map $\hat i_x: \hat W^-_x\to M$, where $\hat W^-_x$ is a compact smooth manifold with corners whose interior is $W^-_x$.
Hence, it follows that  $$\int_{W^-_y} e^{tf}\omega^q_x(t)= \int_{\hat W^-_y} (\hat i_x)^\ast (e^{tf}\omega^q_x(t)).$$
This implies the uniform convergence of the integral (\ref{E1.1}) and the analyticity of $A^q(x,y)(t)$  and of $a^q(M,g,f)(t).$
To conclude the independence on the choices of $\omega^q_\alpha (t)$  it suffices to note that the matrices $|| A^q(\cdots) (t)||$ for two choices differ one from the other by composition by an orthogonal matrix, hence have the same determinant up to sign.
 The results in  \cite{BFK} shows also that  the partition  $M= \bigcup_x W^-_x$ provides a CW structure of $M$ with open cells $W^-_x$
 and if one equips each cell $W^-_x$ with the orientation $O_x$  and one denotes by $C^q$  the $\mathbb R$-vector space  $Maps (C_q(f),\mathbb R)$, then  $Int^q(t): \Omega^q(M) \rightarrow C^q$ defined by
$$Int^q(t)(\omega)(x) = \int _{W^-_x} \omega$$ provides  a linear map, and the collection of the linear maps
$Int^q(t): \Omega^q(M) \rightarrow C^q$ define a map of cochain complexes
$Int^\ast (t): (\Omega^\ast (M), d^\ast(t)) \rightarrow (C^{\ast}, \partial ^\ast)$,  which by de Rham theorem is a quasi- isomorphism. Of course this involves the explicit description of the corner structures of $\hat W^-_x.$
The cochain complex morphism $Int^\ast (t)$ is  clearly analytic in $t$ and restricts to $(\Omega^\ast (M), d^\ast(t))$ an analytic family of quasi-isomorphisms $(\Omega^\ast_{vs}(M), d^\ast (t))\rightarrow (C^\ast, \partial ^\ast).$

If one considers $\omega^q_x(t)$'s as a basis for $\Omega^q_{vs} (M)(t)$ and the characteristic functions of the set $Cr_q(f)$ as a basis for $C^q$, then one realizes that the matrix representation of $Int^q(t)$ is exactly the matrix $||A^q(x,y)(t) ||.$

For $t$ large enough Witten-Helffer-Sj\"ostrand  results imply that $Int^q(t)$ restricted to $\Omega^q_{vs}(M)(t)$ is an isomorphism,
(for details \cite{BFKM} theorem 5.5. item 5, or \cite {B} theorem 3.1), which shows that, for $t$ large enough, $a^q(M,g,f)(t)\ne 0$. In view of analyticity   item 2 follows as stated.
\vskip .1in

Proof of Theorem (\ref{T1.3}):

 First observe the following facts.
\begin{enumerate}

\item For an isomorphism $\varphi: (V, ~ \langle ~, ~ \rangle_V)\to (W, ~ \langle~,~\rangle_W)$ between two finite dimensional vector spaces equipped with scalar product, let
$\Vol (\varphi):= \sqrt{\ddet(\varphi^\sharp \cdot \varphi)^{1/2}}= \sqrt{\ddet(\varphi \cdot \varphi^\sharp)^{1/2}}$ with $\varphi^\sharp$ the adjoint of $\varphi.$
\item
If  $\varphi (t):(V(t), ~\langle~,~\rangle_{V(t)})  \to  (W(t),~ \langle~,~\rangle_{W(t)})$
is a continuous/analytic family of isomorphisms between  finite dimensional vector spaces  equipped with scalar products, \footnote{for example $V(t)$ resp.$W(t)$ appear as images in $\mathcal V$ resp.$\mathcal W,$  of an analytic/continuous family of bounded projectors $P(t): \mathcal V \to \mathcal V$ resp. $Q(t): \mathcal W\to \mathcal W$ for $\mathcal V$ resp  $\mathcal W$ topological vector spaces; this give meaning to "analytic family"}
then the function $\Vol(\varphi (t))$ is continuous/analytic in $t.$
 \item For a cochain complex $\mathcal C= (C^\ast, d^\ast)$ of finite dimensional vector spaces equipped with scalar products
$$ \mathcal C: \  \xymatrix {0\ar[r] &(C^0, \langle~,~\rangle_0)\ar[r]^{d^0} &(C^1, \langle~,~\rangle_1)\ar[r]^{d^1} &\cdots &(C^n, \langle~,~ \rangle_n)\ar[r]&0}, $$
one denotes by $\Delta^q_{\mathcal C}:= \delta^{q+1}\cdot d^q + d^{q-1}\cdot\delta^q,$ \  $\delta $ the adjoint of $d,$ and by $\det' \Delta_{{\mathcal C}}^q\ne 0$  the product of nonzero eigenvalues  of $\Delta_{{\mathcal C}}^q.$
The product
$$T(\mathcal C):= \prod  ({\det}' \Delta_{{\mathcal C}}^q)^{\frac{1}{2} q (-1)^ {q+1}}$$
is referred to as the torsion of $\mathcal C.$ Here $\det \Delta'$ denotes the product of nonzero eigenvalues of $\Delta.$
For a continuous/analytic family of cochain complexes $\mathcal C(t)= (C^\ast(t), d^\ast(t)) $  such that $\dim C^q(t)$ and $\dim H^q(\mathcal C(t))$ are constant in $t$
for any $q$,  the function $T(\mathcal C(t))$ is continuous/analytic in $t.$
\newline The verifications of items (2) and (3)  above are straightforward from definitions.

\item Suppose that $\varphi  : {\mathcal C}_{1} \to {\mathcal C}_{2}$  is a morphism of cochain complexes of finite dimensional vector spaces with scalar products, where $\mathcal C_{i}= (C^\ast_{i}, d^\ast _{i})$, $i=1,2$, and   $\varphi = \{ \varphi^q : C_1^q \to  C_2^q \}.$
Suppose that $\varphi^q$ is an isomorphism for any $q$. Then  $\varphi$ induces
the isomorphism $H^q(\varphi): H^q({\mathcal C}_{1}) \to H^q({\mathcal C}_{2})$ between vector spaces equipped  with induced scalar product.   Let
$$\Vol (\varphi): =\prod (\vol (\varphi^q))^{(-1)^{q}}$$ and $$\Vol (H(\varphi)): = \prod (\vol (H^q(\varphi)))^{(-1)^{q}} .$$ As verified in  \cite {BFK2} Proposition 2.5 one has

\begin{equation}\label{E6.3}
T(\mathcal C_2) / T(\mathcal  C_1)= \Vol(H(\varphi)) / \Vol (\varphi).
\end{equation}
\item  For a continuous/analytic family of isomorphisms $\varphi (t): \mathcal C_1(t)\to \mathcal C_2(t),$\  $t\in \mathbb R,$ with $\dim C^q_1(t)= \dim C^q_2(t)$ and
$\dim H^q({\mathcal C}_1(t))= \dim H^q({\mathcal C}_2 (t))$  constant in $t,$ the real-valued functions
$T({\mathcal C}_1(t))$,  $T({\mathcal C}_2(t))$,  $\Vol (\varphi(t))$, $\Vol (H(\varphi(t)))$ are nonzero and continuous/analytic.
\end{enumerate}

We consider  $\varphi (t)= Int^\ast (t) : (\Omega^\ast_{vs}(M)(t), d^*(t)) \to (C^\ast,\partial ^\ast)$  with $\ast= 0,1, \cdots, \dim M$.  The first cochain complex is equipped with the scalar products defined by the metric $g$, and the second with the unique scalar product which makes the characteristic functions of the critical points orthonormal.
In view of item (4) the function
$$\frac{T(\Omega^\ast_{vs}(M)(t), d^\ast(t))\cdot \Vol (H(\varphi(t))}{T(C^\ast,\partial ^\ast) }$$
is a strictly positive analytic  function and in view of (4) agrees with $a(t)$ for all  $t$ but the finite collection which might be a zero or a  pole for $a(t).$  Hence the  meromorphic function $a(t)$ has no zeros and no poles.
This establishes  Theorem \ref{T1.3} part 1. Together with  (\ref{E6.3}) it also implies
$$\frac{T(\Omega^\ast_{vs}(M)(t), d^\ast(t))}{a(t)} \cdot \Vol (H(\varphi(t))  =  T(C^\ast,\partial^\ast).$$
Evaluation at $t=0$ combined with the observation that $\mathbb Tor(M)= T(C^\ast,\partial^\ast)$    implies
$$\frac{T(\Omega^\ast_{vs}(M), d^\ast)}{a(0)} \cdot  \Vol (H(\varphi(0)))= \mathbb Tor(M).$$
Taking  "$\log$", one derives Theorem \ref{T1.3} part 2.

\section {Poincar\'e duality  for the virtually small spectral package } \label {S3}

For a closed oriented Riemannian manifold  $(M,g)$ and smooth function $f:M\to \mathbb R$
one has the Hodge star operator $\star^q:\Omega^q(M)\to \Omega^{n-q}(M)$ which satisfies the following properties.

\begin{enumerate}
\item $\star^{n-q} \cdot \star^q = (-1)^{q(n-q)} \Id,$
\item $(-1)^{q(n-q)}\star^q  \Delta^q \ \star^{n-q} = \Delta^{n-q},$
\item $(-1)^{q(n-q)}\star^q  \Delta_{g,f}^q(t)\ \star^{n-q} =\Delta^{n-q}_{g, -f}(t),$
\item $\Delta^q_{g,f}(-t)= \Delta^q_{g, -f} (t).$
\end{enumerate}

As a consequence, the Hodge operator $\star^q$ identifies the $q$-spectral package
of $~(M, g, f)~$  with the $(n-q)$-spectral package of $(M, g, -f)$ and the
$q$-virtually small spectral package

\begin{eqnarray*}
\left\{ \lambda^q_{g, f, x}(t), ~ \omega^q_{g, f, x}(t), ~ x\in Cr_q(f) \right\}
\end{eqnarray*}

\noindent
of $(M, g, f)$ with the $(n-q)$-virtually small spectral package

\begin{eqnarray*}
\left\{\lambda^{n-q}_{g, -f, x}(t), ~ \omega^{n-q}_{g, -f, x}(t), ~ x\in Cr_{n-q}(-f) \right\}
\end{eqnarray*}

\noindent
of $(M, g, -f)$.
More precisely,  it holds that for $x\in Cr_{q} (f)= Cr_{n-q} (-f)$,

\begin{eqnarray*}
\lambda^{n-q}_{g,-f,x}(t) ~ = ~ \lambda^{q}_{g,f,x}(t), \qquad
\omega^{n-q}_{g,-f, x}(t) ~ = ~ \star^{q}  \omega^{q}_{g,f,x}(t).
\end{eqnarray*}

\vskip .1in

{\bf Poincar\'e duality}

The above  identification can be regarded as a refinement of the Poincar\'e duality  which states that  $\beta^q(M)$ viewed as the multiplicity of the eigenvalue $0$ of $\Delta^q_{g,f}(0)$
is equal to $\beta^{n-q}(M)$ viewed as the multiplicity of the eigenvalue $0$ of $\Delta^{n-q}_{g,-f}(0).$
\vskip .1in

 \section {Virtually small eigenvalues versus the smallest eigenvalues} \label {S4}

We are going to show that the virtually small eigenvalues may not be equal to the smallest eigenvalues by giving an example.
For $\mathbb S^1:= \mathbb R/ (2\pi \mathbb Z)$, we consider a torus $M= \mathbb S^1\times \mathbb S^1$  equipped with the flat metric $g_0$ induced from the canonical metric on $\mathbb R\times \mathbb R$.
The function
 $f(\theta_1, \theta_2)= \sin(2\theta_1) +\sin(2\theta_2)$,  $(\theta_1, \theta_2) \in \mathbb S^1\times \mathbb S^1$,  is a Morse function on $M$ having
four critical points of index 0, four critical points of index $2$ and eight critical points of  index  $1$.
The sequence of eigenvalues of  $\Delta^0_{g_0}$ in increasing order is $0\leq 1 \leq 1 \leq 1\leq  1\leq 4\cdots$.
The virtually small $0$-eigenvalues consist of four real numbers. The first one is $\lambda_1(0) = 0$, the  next two are $\lambda_2(0)= \lambda_3(0)  \geq 1 $
and the final one is $\lambda_4(0)= 2 \mu_2(0)$, where $\mu_{2}(0) \geq 1$.
This shows that the virtually small eigenvalues are not the same as the smallest eigenvalues.

\begin {proof}
We first observe that $h:\mathbb S^1\to \mathbb R$ given by $h(\theta)= \sin 2\theta$ is a Morse function on ${\mathbb S}^{1}$, whose Witten Laplacian $\Delta^0(t)$ is
$$\Delta^0(t)= -\partial ^2/\partial \theta^2 +4t \sin 2\theta  + 4t^2 (\cos 2\theta)^2.$$

\noindent
Then, $\Delta^0(t)$ has two virtually small eigenvalue branches. One of them is $\mu_1 (t)\equiv 0$
and the other is $\mu_2(t) >0$ because $\beta^0(\mathbb S^1)=1$ and $h$ has two critical points of index $0.$
Observe that the eigenvalues of $\Delta^0(0)$ are

$$0, ~ 1, ~ 2^2, ~ 3^2, ~ \cdots n^2, ~ \cdots,$$

\noindent
where $0$ has the multiplicity $1$ and all others have multiplicity $2$.
Since $\mu_1(0)$ and  $\mu_2(0)$ are among the above eigenvalues, one has $\mu_2(0)\geq 1.$
In view of the definition of $f (\theta_1, \theta_2)$  and that of $\Delta^0(t)$ for $(\mathbb S^1\times \mathbb S^1, g_0, f)$,  the four virtually small eigenvalues branches $\lambda_1(t),$ $ \lambda_2(t), \lambda_3(t) ,\lambda_4(t)$ are of the form $\mu_i (t)+ \mu_j(t)$ with $i,j \in \{1,2\}$, and hence the virtually small $0$-eigenvalues are $\lambda_1(0)=0, \lambda_2(0)= \lambda_3(0)= \mu_1(0)$ and $\lambda_4(0)= 2 \mu_1(0).$
 \end {proof}

When $M^{2}$ is a $2$-dimensional oriented closed Riemannian manifold,
the above example shows that the  virtually small spectral package of $\Delta^{0}$ is not the same as the collection of the smallest eigenvalues.
By the Poincar\'e duality, the nonzero eigenvalues of $\Delta^1$ on $M^{2}$ are two times of the nonzero eigenvalues of $\Delta^0$ or $\Delta^2$, and hence
the result remains  the same for $\Delta^1.$

\vskip .2in
{\small
Dan Burghelea

Department of Mathematics,

Ohio State University,
Columbus, OH 43210, USA}

\vskip .2in
{\small
Yoonweon Lee

Department of Mathematics Education,

Inha University,
Incheon, 22212, Korea}

 \end{document}